\documentclass[12pt,reqno,a4wide]{amsart}

\usepackage{tikz} 
\usetikzlibrary{positioning}
\usepackage{calrsfs}

\allowdisplaybreaks

\title{Enumeration of Deutsch paths by  the adding-a-new-slice method and applications}

\author{Helmut Prodinger}
\address{Helmut Prodinger\\
	Department of Mathematical Sciences\\
	Stellenbosch University\\
	7602 Stellenbosch\\
	South Africa}
\email{hproding@sun.ac.za}
\begin{document}
\begin{abstract}
A variation of Dyck paths allows for down-steps of arbitrary length, not just one. This is motivated by ideas due to Emeric Deutsch.  We use the adding-a-new-slice technique and the kernel method to compute the number of maximal runs of up-step runs of length 1
and a subclass of Deutsch paths satisfying a condition that was stipulated by R. Stanley for Dyck paths. 
	\end{abstract}

\maketitle

\section{Introduction}

A Dyck path consists of up-steps $(1,1)$ and down-steps $(1,-1)$. They start at the origin, end at the origin and never go below the $x$-axis.
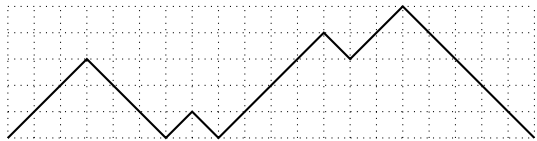
\begin{center}

\begin{tikzpicture}[scale=0.35]
	
	\draw[step=1.cm,black,dotted] (-0.0,-0.0) grid (20.0,5.0);
	\draw[thick] (0,0) to (1,1) to (2,2) to (3,3) to (4,2) to (5,1) to (6,0)to (7,1)to (8,0)to (9,1)to (10,2)to (11,3)to (12,4)to (13,3)to (14,4)to (15,5)to (16,4)to (17,3)to (18,2)to (19,1)to (20,0);
\end{tikzpicture}
\end{center}
The path in the example is a Dyck path of length 20, and has 4 `mountains', i.e. a maximal sequence of up-steps followed by a maximal sequence of down-steps.
Sometimes we also allow a path to end at a level different from zero. To emphasize the fact that the end-level is unspecified, we use the word \emph{open}.
In \cite{BF}, this is called a \emph{meander}, whereas the version that returns to the $x$-axis is called an \emph{excursion}.

Deutsch paths are like Dyck paths, but extra down-steps of the form $(1,-j)$, for any $j\ge2$, are also allowed.
They were analyzed recently in \cite{Prodinger-Deutsch}. Here, we want to enumerate them in a different manner
which is quite versatile when certain parameters of Deutsch paths should be analyzed.

We decompose a Deutsch path into maximal runs of up- resp.\ down-steps. First, we restrict our attention to the instance when the path ends with down-steps. If a path is closed (=ends on the $x$-axis), this happens anyway (except for the empty path), but for open-ended paths, the last step of the path might be an up-step.

The technique we are using can be found in \cite{FP}.

We consider two applications: 
It is well-known that the $n$-th Catalan number $C_n =\frac{1}{n+1}\binom{2n}{n}$
enumerates Dyck paths of
length $2n$. In \cite{Stanley}, Stanley lists a variety of other combinatorial interpretations of the Catalan
numbers, one of them being the number of Dyck paths from $(0, 0)$ to $(2n+2, 0)$ such that any
maximal sequence of consecutive $(1,-1)$ steps ending on the $x$-axis has odd length. At this
point it is interesting to note that there are more subclasses of Dyck paths, also enumerated
by Catalan numbers, that are defined via parity restrictions on the length of the returns to the
$x$-axis as well (see, e.g., \cite{Stanley}). This restriction of Dyck paths that leads again to Catalan numbers
was further investigated in  \cite{HP}. In this paper, we consider Deutsch paths with the property
that each maximal run of down-steps to the $x$-axis starts at an odd level. Unfortunately, in this context,
this does not lead to any known or  nice numbers.

In the last section, we count the number of maximal sequences of up-steps consisting of only one up-step in Deutsch paths.

In the paper \cite{Estland}, the authors counted the number of maximal runs of up-steps of length one in Dyck paths. This was greatly extended in \cite{HHP};
the method of choice in the first incarnation of this paper was indeed the adding-a-new-slice technique combined with the kernel method.  We keep our analysis in the instance of Deutsch paths  simple and only consider the basic case, leaving more general considerations for later.

\section{Enumeration of Deutsch paths by the adding-a-new-slice technique}

We consider a generating function $F_k(z,u)$, where $k$ is the number of `mountains' (runs of up-steps, followed by runs of down-steps), $z$ marks the length of the path and $u$ is used to remember the last level reached (coefficient of $u^j$). To be more precise, the coefficient of $z^nu^j$ in $F_k(z,u)$, for which we write
$[z^nu^j]F_k(z,u)$ is the number of Deutsch paths with $k$ mountains, length $n$, and ending on level $j$.

The extra down-steps require some preparations. If one want to go down by $h\ge1$ levels, and do this in $n\ge1$ steps, it can be done in
\begin{equation*}
[z^h]\Bigl(\frac{z}{1-z}\Bigr)^n=[z^{h-n}](1-z)^{-n}=\binom{h-1}{h-n}
\end{equation*}
ways. (The notation of $[z^h]f(z)$ refers, as just said, to the coefficient of $z^n$ in a generating function $f(z)$.)
To understand this, we are looking at the coefficient of $z^h$ in
\begin{equation*}
\underbrace{(z+z^2+\cdots)(z+z^2+\cdots)\dots(z+z^2+\cdots)}_{n \text{ copies}},
\end{equation*}
and any choice of $i_1+\cdots+i_n=h$ is a possibility to go down by $h$ levels in $n$ steps.
 If one wants to keep $n$ variable, using a generating function, we get, using the binomial theorem,
\begin{equation}\label{down-h}
\sum_{n=1}^{h}\binom{h-1}{h-n}z^n=z(1+z)^{h-1}.
\end{equation}
This will also be used later. And now, for modeling a new mountain,  we have to compute this:
\begin{equation*}
\sum_{k>i}z^{k-i}\sum_{0\le j<k}u^jz(1+z)^{k-1-j}=\frac{z^2(1+z)^{i+1}}{(1-z-z^2)(1+z-u)}-\frac{z^2u^{i+1}}{(1-zu)(1+z-u)}.
\end{equation*}
Consequently we have
\begin{equation*}
F_{k+1}(z,u)=\frac{z^2(1+z)}{(1-z-z^2)(1+z-u)}F_k(z,1+z)-\frac{z^2u}{(1-zu)(1+z-u)}F_k(z,u).
\end{equation*}
With 
\begin{equation*}
\Phi(z,u):=\sum_{k\ge0}F_k(z,u),
\end{equation*}
(arbitrary number of mountains) this leads to
\begin{equation*}
	\Phi(z,u)-1=\frac{z^2(1+z)}{(1-z-z^2)(1+z-u)}\Phi(z,1+z)-\frac{z^2u}{(1-zu)(1+z-u)}\Phi(z,u)
\end{equation*}
or
\begin{equation*}
	\frac{z(u-u_1)(u-u_2)}{(1-zu)(1+z-u)}\Phi(z,u)=1+\frac{z^2(1+z)}{(1-z-z^2)(1+z-u)}\Phi(z,1+z),
\end{equation*}
with
\begin{equation*}
	u_{1,2}=\frac{1+z\pm\sqrt{1-2z-3z^2}}{2z}.
\end{equation*}
A naive approach to solve this functional equation would be to plug in $u=1+z$ and solve. However, this leads to a void equation. The kernel method is a way to still use
the idea of plugging in, but first a `bad' factor has to be removed from numerator and denominator. One of two factors $(u-u_1(z))$, $(u-u_2(z))$ is bad, since, when appearing in the denominator, does not lead to a valid power series expansion around $(0,0)$. In this instance, $u-u_2(z)\sim u-z+\cdots$ and this is the bad factor. Rephrasing, setting
$u:=u_2(z)$ leads to an equation, since the bad factor must be cancelled from both, numerator and denominator.

As it was shown already in \cite{Prodinger-Deutsch}, we are in the Motzkin world.
Recall that the generating function of Motzkin-paths (Dyck paths, but horizontal steps are also allowed) is given by
\begin{equation*}
M(z)=\frac{1-z-\sqrt{1-2z-3^2}}{2z^2}.
\end{equation*}
Now we use arguments from the kernel method \cite{kernel}, as just recapitulated.
The factor $(u-u_2)$ must also be a factor of the right-hand side, otherwise there would not be a power series expansion around
$z=0$. This leads to
\begin{equation*}
\frac{1+\frac{z^2(1+z)}{(1-z-z^2)(1+z-u)}\Phi(z,1+z)}{u-u_2}=-\frac1{1+z-u}.
	\end{equation*}
	Further simplification leads to
\begin{equation*}
	\Phi(z,u)=-\frac{(1-zu)}{z(u-u_1)}=\frac{1-zu}{zu_1(1-u/u_1)}.
\end{equation*}
	Setting $u=0$ means that the path ends on the $x$-axis:
\begin{equation*}
	\Phi(z,0)=\frac{1}{zu_1}=\frac{1+z-\sqrt{1-2z-3z^2}}{2z(1+z)}=
	1+{z}^{2}+{z}^{3}+3{z}^{4}+6{z}^{5}+15{z}^{6}+36{z}^{7}+\cdots.
\end{equation*}
The coefficients of this series are sometimes called \emph{Riordan numbers}, see A005043 in \cite{OEIS}.
As an illustration, we provide the 6 Deutsch paths (returning to the $x$-axis) of length 5:
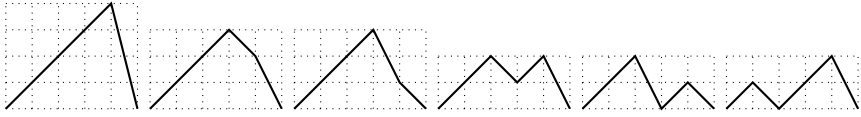
\begin{center}
	\begin{tikzpicture}[scale=0.35]
		\draw[step=1.cm,black,dotted] (-0.0,-0.0) grid (5.0,4.0);
		%\draw[step=1.cm,black,dotted] (-0.0,-0.0) grid (9.0,3.0);
		\draw[thick] (0,0) to (1,1) to (2,2) to (3,3) to (4,4) to (5,0);
	\end{tikzpicture}
\begin{tikzpicture}[scale=0.35]
	\draw[step=1.cm,black,dotted] (-0.0,-0.0) grid (5.0,3.0);
	%\draw[step=1.cm,black,dotted] (-0.0,-0.0) grid (9.0,3.0);
	\draw[thick] (0,0) to (1,1) to (2,2) to (3,3) to (4,2) to (5,0);
\end{tikzpicture}
	\begin{tikzpicture}[scale=0.35]
		\draw[step=1.cm,black,dotted] (-0.0,-0.0) grid (5.0,3.0);
		%\draw[step=1.cm,black,dotted] (-0.0,-0.0) grid (9.0,3.0);
		\draw[thick] (0,0) to (1,1) to (2,2) to (3,3) to (4,1) to (5,0);
	\end{tikzpicture}
\begin{tikzpicture}[scale=0.35]
	\draw[step=1.cm,black,dotted] (-0.0,-0.0) grid (5.0,2.0);
	%\draw[step=1.cm,black,dotted] (-0.0,-0.0) grid (9.0,3.0);
	\draw[thick] (0,0) to (1,1) to (2,2) to (3,1) to (4,2) to (5,0);
\end{tikzpicture}
\begin{tikzpicture}[scale=0.35]
	\draw[step=1.cm,black,dotted] (-0.0,-0.0) grid (5.0,2.0);
	%\draw[step=1.cm,black,dotted] (-0.0,-0.0) grid (9.0,3.0);
	\draw[thick] (0,0) to (1,1) to (2,2) to (3,0) to (4,1) to (5,0);
\end{tikzpicture}
\begin{tikzpicture}[scale=0.35]
	\draw[step=1.cm,black,dotted] (-0.0,-0.0) grid (5.0,2.0);
	%\draw[step=1.cm,black,dotted] (-0.0,-0.0) grid (9.0,3.0);
	\draw[thick] (0,0) to (1,1) to (2,0) to (3,1) to (4,2) to (5,0);
\end{tikzpicture}

\end{center}

\medskip

	As seen  already in \cite{Prodinger-Deutsch}, the substitution $z=\dfrac{v}{1+v+v^2}$ makes everything prettier:
\begin{equation*}
	\Phi(z,0)=\frac{1+v+v^2}{1+v}.
\end{equation*}
The function
\begin{equation*}
\Phi(z,u)\frac1{1-zu}
\end{equation*}
describes Deutsch paths that can also end with up-steps. And if one replaces now $u:=1$, we get so-called open Deutsch paths, that can end at any level:
\begin{equation*}
	\Phi(z,1)\frac1{1-z}=\frac{1}{z(u_1-1)}=1+v+v^2,
\end{equation*}
which also enumerates Motzkin paths. This was explained via a bijection in \cite{Prodinger-Deutsch}. Let us emphasize that $	\Phi(z,1)$ is
the generating function of open Deutsch paths ending with down-steps, and $\frac1{1-z}$ is the generating function of a possible sequence of up-steps at the end.

\section{Deutsch paths satisfying a  condition by Stanley}
As mentioned above, the function
\begin{equation*}
\Phi(z,u)\frac1{1-zu}
\end{equation*}
describes Deutsch paths that can also end with up-steps. Consequently,  
\begin{equation*}
\mathcal{G}(z,u)=	\Phi(z,u)\frac{zu}{1-zu}=\frac{u}{u_1(1-u/u_1)}=\sum_{k\ge1}\frac{u^k}{u_1^k}
\end{equation*}
is the generating function of paths (`good' paths) ending with an up-step. From this we see that the good paths ending on level $k$ have
generating function $1/u_1^k$.

In the spirit of Stanley, we now compute good Deutsch paths, ending on the odd level $2k+1$, and return after that for the \emph{first time} in a series of down-steps to the $x$-axis:
\begin{equation*}
z\sum_{k\ge0}\frac{1}{u_1^{2k}}\cdot z(1+z)^{2k} =\frac{z^2}{1-(1+z)^2/u_1^2}.
\end{equation*}
To clarify, the first $z$ in this formula is responsible for a first extra up-step, making sure that returns to the $x$-axis occur only at the end.
The other factor $z$ belongs to $z(1+z)^{2k}$ and originates from formula \eqref{down-h}.

The final step is to consider an arbitrary sequence of such paths, viz.
\begin{equation*}
\dfrac{1}{1-\frac{z^2}{1-(1+z)^2/u_1^2}}=\frac{3+z-\sqrt{1-2z-3z^2}}{2(1+z)}=\frac{1+2v+2v^2}{(1+v)^2}.
\end{equation*}
This series is
\begin{equation*}
1+{z}^{2}+2 {z}^{4}+2 {z}^{5}+7 {z}^{6}+14 {z}^{7}+37 {z}^{8}+90
 {z}^{9}+233 {z}^{10}+\cdots,
\end{equation*}
and the coefficients do not bare any significance to Motzkin numbers.

For interest, here are the $7$ objects of length $6$, satisfying the Stanley-condition that each down-run to the $x$-axis starts at an odd level.

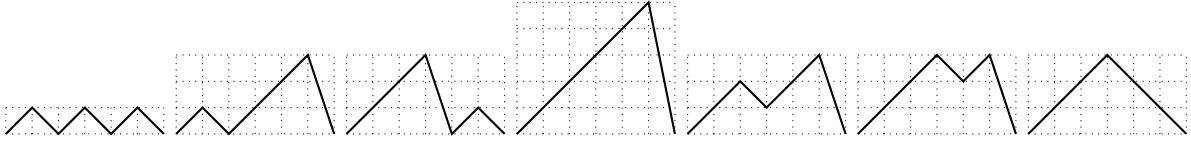
\begin{center}
\begin{tikzpicture}[scale=0.35]
\draw[step=1.cm,black,dotted] (-0.0,-0.0) grid (6.0,1.0);
%\draw[step=1.cm,black,dotted] (-0.0,-0.0) grid (9.0,3.0);
\draw[thick] (0,0) to (1,1) to (2,0) to (3,1) to (4,0) to (5,1) to (6,0);
\end{tikzpicture}
\begin{tikzpicture}[scale=0.35]
\draw[step=1.cm,black,dotted] (-0.0,-0.0) grid (6.0,3.0);
%\draw[step=1.cm,black,dotted] (-0.0,-0.0) grid (9.0,3.0);
\draw[thick] (0,0) to (1,1) to (2,0) to (3,1) to (4,2) to (5,3) to (6,0);
\end{tikzpicture}
\begin{tikzpicture}[scale=0.35]
\draw[step=1.cm,black,dotted] (-0.0,-0.0) grid (6.0,3.0);
%\draw[step=1.cm,black,dotted] (-0.0,-0.0) grid (9.0,3.0);
\draw[thick] (0,0) to (1,1) to (2,2) to (3,3) to (4,0) to (5,1) to (6,0);
\end{tikzpicture}
\begin{tikzpicture}[scale=0.35]

%\draw[step=1.cm,black,dotted] (-0.0,-0.0) grid (9.0,3.0);
\draw[step=1.cm,black,dotted] (-0.0,-0.0) grid (6.0,5.0);
\draw[thick] (0,0) to (1,1) to (2,2) to (3,3) to (4,4) to (5,5) to (6,0);
\end{tikzpicture}
\begin{tikzpicture}[scale=0.35]

\draw[step=1.cm,black,dotted] (-0.0,-0.0) grid (6.0,3.0);
\draw[thick] (0,0) to (1,1) to (2,2) to (3,1) to (4,2) to (5,3) to (6,0);
\end{tikzpicture}
\begin{tikzpicture}[scale=0.35]

\draw[step=1.cm,black,dotted] (-0.0,-0.0) grid (6.0,3.0);
\draw[thick] (0,0) to (1,1) to (2,2) to (3,3) to (4,2) to(5,3) to (6,0);
\end{tikzpicture}
\begin{tikzpicture}[scale=0.35]

\draw[step=1.cm,black,dotted] (-0.0,-0.0) grid (6.0,3.0);
\draw[thick] (0,0) to (1,1) to (2,2) to (3,3) to (4,2) to (5,1) to (6,0);
\end{tikzpicture}

\end{center}
The first and the last   are Dyck-paths, and 2 is indeed the Catalan number $C_2$.

	\section{Counting runs of single up-steps}

The approach is quite similar to the previous sections; however, we use a third variable, $t$, to count the up-runs of length one.

We have to compute this (for $t=1$, it coincides with the previous computation):
\begin{equation*}
	\sum_{k>i+1}z^{k-i}\sum_{0\le j<k}u^jz(1+z)^{k-1-j}+tz\sum_{0\le j\le i}u^jz(1+z)^{i-j},
\end{equation*}
which leads to
	\begin{equation*}
		F_{k+1}(z,u)=\alpha F_k(z,u)+\beta F_k(z,1+z)
	\end{equation*}
	with
	\begin{equation*}
\alpha={\dfrac {{z}^{2}u \left( -zu+tzu-t \right) }{ ( 1-zu ) 			( 1+z-u ) }}
	\quad\text{and}\quad\beta=-{\dfrac {{z}^{2} ( 1+z )  ( -z-{z}^{2}+tz+t{z}^{2}-t) }{( 1-z-{z}^{2} )  ( 1+z-u )  }}.
	\end{equation*}	
	This leads to
	\begin{equation*}
		\Phi(z,u)-1=\alpha\Phi(z,u)+\beta\Phi(z,1+z)
	\end{equation*}
	or
	\begin{equation*}
		\frac{z(1+z^2-tz^2)(u-u_1)(u-u_2)}{(1-zu)(1+z-u)}\Phi(z,u)=1+\beta\Phi(z,1+z),
	\end{equation*}
	
	The two roots are now
	\begin{equation*}
		u_{1,2}={\frac {-t{z}^{2}+z+1+{z}^{2}\pm\sqrt {{t}^{2}{z}^{4}+2 t{z}^{3}-
					2 t{z}^{2}+2 t{z}^{4}-{z}^{2}-2 z-2 {z}^{3}+1-3 {z}^{4}}}{2z
				( 1+{z}^{2}-t{z}^{2}) }}.
	\end{equation*}
	Simplification, after dividing out the factor $u-u_2$ from the equation, leads to
\begin{equation*}
\Phi(z,u)=\frac{(1-zu)}{z(1+z^2-tz^2)(u_1-u)}.
\end{equation*}
This time, we confine ourselves to the instance $u=0$, i.e., Deutsch paths returning to the $x$-axis. We get
\begin{align*}
	\Phi(z,&0)=\frac{1}{z(1+z^2-tz^2)u_1}\\&
={\frac {-t{z}^{2}+z+1+{z}^{2}-\sqrt {{t}^{2}{z}^{4}+2 t{z}^{3}-2
			 t{z}^{2}+2 t{z}^{4}-{z}^{2}-2 z-2 {z}^{3}+1-3 {z}^{4}}}{2z
		( 1+z )  ( {z}^{2}-t{z}^{2}+1 ) }}\\&
=1+t{z}^{2}+{z}^{3}+ ( {t}^{2}+2 ) {z}^{4}+ ( 3+3 t) {z}^{5}+ ( 7+7 t+{t}^{3} ) {z}^{6}+ 
( 17+13 t+6 {t}^{2} ) {z}^{7}+\cdots.
\end{align*}
Once one has this generating function, one can state many results as a corollary. We will only provide one such result, namely,  the average of the parameter labelled by the variable $t$. So, we differentiate $\Phi(z,0)$ w.r.t. $t$ and then set $t:=1$.
This leads to
\begin{equation*}
\frac{v^2}{(1-v)(1+v)^2(1+v+v^2)}.
\end{equation*}
We prefer the factored form since the term $1+v+v^2$ immediately reminds us of the Motzkin-world.

One could even read off the coefficients from this, but this would lead to a sum, so we refrain from doing this. However, we are interested in asymptotics. We will use singularity analysis, as is now customary, see \cite{FS}. In particular Example VI.3. on page 396 discusses the basics of the asymptotics of Motzkin numbers.

The relevant singularity is at $z=\frac13$, and we find, as $z\to\frac13$ (the square root vanishes at this value, and it is the closest value to the origin with this property).
Since 
\begin{equation*}
v=\frac{1-z-\sqrt{1-2z-3z^2}}{2z},
\end{equation*}
the singularity translates into $v=1$, and one can compute a local expansion:
\begin{equation*}
v\sim 1-\sqrt{3}\sqrt{1-3z} \quad \Longleftrightarrow\quad 1-v \sim -\sqrt{3}\sqrt{1-3z}.
\end{equation*}

Furthermore,
\begin{equation*}
\frac{v^2}{(1-v)(1+v)^2(1+v+v^2)}\sim\frac{1}{12(1-v)}  \sim \frac{\sqrt3}{36\sqrt{1-3z}}.
\end{equation*}
Therefore
\begin{equation*}
[z^n]	\frac{v^2}{(1-v)(1+v)^2(1+v+v^2)}\sim [z^n]\frac{\sqrt3}{36\sqrt{1-3z}}\sim\frac{\sqrt3}{36}3^n\frac1{\sqrt{\pi n}}.
	\end{equation*}
The table on page 388 in \cite{FS} is particularly useful to translate from local expansions to asymptotics of coefficients.

This needs to be divided by the total number of such paths (the coefficient of $z^n$ in $\Phi(z,0)$), viz.
\begin{equation*}
[z^n]\frac{1+v+v^2}{1+v}\sim[z^n]\Big(\frac32-\frac{3\sqrt3}{4}\sqrt{1-3z}\Big)\sim \frac{3\sqrt3}{8}3^n\frac1{\sqrt{\pi} n^{3/2}}.
	\end{equation*}
The quotient is
\begin{equation*}
\sim\frac{2n}{27}=0.074 n.
\end{equation*}
So, a Deutsch path of length $n$ has about $0.074 n$ up-runs of length 1. Many such results could be derived with some patience and a computer.

\section{Concluding remarks}

We have demonstrated the usefulness of the adding-a-new-slice technique \cite{FP} in the instance of the relatively new class of Deutsch paths, 
introduced in \cite{Prodinger-Deutsch}. Whenever it is possible to describe how a new mountain is created, a functional equation can be written. 
A variation of the idea of Deutsch paths is to use down-steps $-1,-3,-5,-7,\dots$  instead of $-1,-2,-3,-4,\dots$, as done here.
This situation was analyzed in \cite{Prodinger-special}, but not with the adding-a-new-slice method. This is probably possible, but instead of a quadratic
equation with 2 factors, a cubic equation with 3 factors will play a role, and either one or two of them (depending on the model) need to be cancelled out.

\medspace

\textsc{Acknowledgment.} The insightful comments of one referee resulted in a much more elaborate and hopefully clearer exposition of this material.

%\clearpage

\end{document}